\documentclass{article}[12pt]
\usepackage[margin=1in]{geometry}
\usepackage{amsmath,amsfonts,amsthm}

\newtheorem{assumption}{Assumption}
\newtheorem{theorem}{Theorem}
\newtheorem{lemma}{Lemma}
\newtheorem{proposition}{Proposition}
\newtheorem{remark}{Remark}

\newcommand{\real}{{\mathbb{R}}}
\newcommand{\natch}{{\mathbb{N}}}

\newcommand{\one}{{\mathbb{I}}}
\newcommand{\refeq}[1]{{(\ref{#1})}}
\newcommand{\stochle}{\stackrel{\text{\makebox[0pt]{st}}}{\le}}
\newcommand{\cF}{{\mathcal{F}}}
\newcommand{\cG}{{\mathcal{G}}}
\newcommand{\cN}{{\mathcal{N}}}

\def\Halmos{\mbox{\quad$\square$}}

\begin{document}

\title{Probabilistic bisection converges almost as quickly as stochastic approximation}
\author{Peter I.\ Frazier, Shane G. Henderson, Rolf Waeber}
\date{}
\maketitle


\abstract{%
The probabilistic bisection algorithm (PBA) solves a class of
stochastic root-finding problems in one dimension by successively updating a
prior belief on the location of the root based on noisy responses to 
queries at chosen points. The responses indicate
the direction of the root from the queried point, and are incorrect
with a fixed probability. The fixed-probability assumption is
problematic in applications, and so we extend the PBA to apply when
this assumption is relaxed. The extension involves the use of a
power-one test at each queried
point. We explore the convergence behavior of the extended PBA, showing that it
converges at a rate arbitrarily close to, but slower than, the canonical ``square root''
rate of stochastic approximation.
}%


\section{Introduction.}\label{sec:intro}
The probabilistic bisection algorithm (PBA) \cite{Ho63} is an algorithm
for solving certain one-dimensional root-finding problems where a
query at a single point returns a stochastic response. The stochastic
response is assumed to indicate the direction of the root, and is assumed
correct with some fixed probability $p > 1/2$. The PBA begins with a
prior belief distribution on the location of the root, and repeatedly
queries at the median of the belief distribution, updating the
belief distribution using Bayes' rule. The sequence of medians of the
belief distribution are now known to converge to the true root at an
exponential rate \cite{waefrahen13}. 

The PBA itself, and algorithms based on the PBA, have been applied in contexts including 
target localization \cite{tsiligkaridis2014collaborative},
scanning electron microscopy \cite{SznitmanLucchiFrazier2013},
topography estimation \cite{CaNo08a}, 
edge detection \cite{GoLe03},
and value function approximation for optimal stopping problems \cite{rodriguez2015information}.
A modification of the PBA to the discrete setting was considered in \cite{BuZi74,BuZi75}, as part of a larger body of work in computer science on unknown-item two-person question and answer games with lies \cite{Pe89}.
The PBA is known to be optimal in an average-case sense for an
entropy loss function \cite{JedynakFrazierSznitman2012}.

The PBA is intuitive and the posterior belief distribution provides a
great deal of insight on the location of the root. However, the
assumption that the correctness probability $p$ does not depend on the
location at which a query is performed limits its applicability. To
see why, consider finding the root $x^*$ of a decreasing function $g:[0,1] \to \real$,
where we observe $g(x) + \epsilon$ when querying $g$ at $x \in [0,
1]$, where $\epsilon$ is a normal random variable with mean 0 and
constant variance 1 that is independent of all previous queries and
responses. The query indicates that $x^* > x$ if $g(x) + \epsilon > 0$
and that $x^* < x$ if $g(x) + \epsilon < 0$. If $g(x) \to 0$ as $x \to
x^*$, then the probability that the query result is correct converges
to $1/2$ as $x \to x^*$, violating the assumption of constant
probability of correctness.

We modify the PBA to allow for non-constant
probabilities of correctness. The modification performs
repeated queries at a single point $x$ until the probability of
correctness exceeds a user-selected threshold $p_c > 1/2$. This
can be achieved through the use of tests of power one. Once
the test concludes at $x$, we update the belief distribution
in the same manner as in the PBA, using success probability $p_c$. The
resulting belief distribution is not a true posterior, and so the
existing analysis of the PBA is not valid. Nevertheless, the belief
distribution provides an intuitive representation of our belief about
the location of the root, and can be leveraged to obtain confidence
intervals on the location of the root \cite{Waeber2013}. These
confidence intervals are valid irrespective of how many queries are
performed, i.e., they are not asymptotic, and the resulting confidence
interval widths shrink at an exponential rate in terms of the number
of points at which the function is queried \cite{Waeber2013}.

In this paper, we show that when applying the modified PBA (henceforth
referred to as the PBA for simplicity) on
one-dimensional stochastic root-finding problems, the distance $|X_n - x^*|$
between the median $X_n$ of the belief distribution after running the
test of power one at $n$ points and the true root $x^*$ converges at
an exponential (in $n$) rate to 0. (This result is not implied by the confidence
interval result in \cite{Waeber2013}.) This does not imply that we get
an exponential rate of convergence in {\em wall-clock} time to the root,
however, since the tests of power one have runlengths that increase
without bound as $X_n$ approaches $x^*$. We further show that for any
$\epsilon > 0$, an estimator of
the root based on the first $n$ tests that collectively require $T_n$
total queries (which we take as a reasonable surrogate for wall-clock
time) has an error that is at most of stochastic order
$T_n^{-1/2+\epsilon}$.

By way of contrast, the standard method for stochastic
root-finding problems is stochastic approximation
\cite{RoMo51,KuYi03,PaKi11}. Under mild regularity conditions,
stochastic approximation applied in our setting would have an error of
stochastic order $c^{-1/2}$ based on a wall-clock running time of
$c$. The $\epsilon$ difference in the exponent is striking,
particularly since $\epsilon > 0$ can be arbitrarily small, and
explains the title of this paper. Our analysis does not provide a {\em
  tight} bound on the true stochastic order of the error in PBA, but
we conjecture that it is $T_n^{-1/2} \ln(T_n)$, which is slightly
slower than stochastic approximation.

Even though PBA appears to have a slower asymptotic convergence rate
than stochastic approximation, we still believe the algorithm has
value for several reasons. First, the belief distribution provides
insight on the location of the root. Second, this asymptotic result
does not shed light on how the algorithms compare over wall-clock
running times that are representative in practice. Third, stochastic
approximation requires the selection of a large number of tuning
parameters, notably the ``gain sequence,'' with poor selections
leading to reductions in the rate of convergence; see, e.g.,
\cite[p.~253]{asmgly07}. Fourth, the PBA appears to be readily adapted for use
in parallel computing environments \cite{palfrahen14}, while the
situation is more complicated for stochastic approximation
\cite{kusyin87a,kusyin87b}.

We view the key contributions of this paper to be the extension of the
PBA to a broad class of root-finding problems, and the analysis of the
asymptotic error, both as a function of the number of power-one tests
completed, and as a function of $T_n$, the surrogate for wall-clock
time. The convergence analysis herein is completely different from
that in \cite{waefrahen13} for the original PBA due to necessity, and is
a considerable outgrowth of the original source for these ideas in
\cite{Waeber2013} in which the exponential convergence rate was
conjectured but unproven. It is also completely different from the
analysis in \cite{BuZi74} where it is assumed that the root lies
amongst a finite set of alternatives.

The question of whether the PBA can be successfully applied in
root-finding problems of dimension greater than one remains open. We
do not consider that question here.

The remainder of this paper is organized as follows. The PBA is
reviewed in \S\ref{sec:pba}, and in \S\ref{sec:algorithm} we provide
our extension of the PBA that uses tests of power one to advance at
each step. We establish the exponential (in $n$) rate of convergence of the
sequence of medians in \S\ref{sec:exponential}, and consider
wall-clock time behavior in \S\ref{sec:canonical}. Proofs that are
omitted to this point are collected in \S\ref{sec:proofs}.

\section{Probabilistic Bisection.}\label{sec:pba}
\begin{assumption}
Let $g:[0, 1] \to \real$ and suppose that there exists $x^* \in [0,
1]$ such that $g(x) > 0$ for $x < x^*$ and $g(x) < 0$ for $x >
x^*$.
\end{assumption}
Our goal is to estimate $x^*$, but we make no assumption about
the value of $g(x^*)$.

In this section we assume that we can query the value of $g$ at any
point $x \in (0, 1)$, obtaining a signal $Z(x) \in \{-1, 1\}$, independently of all
previous queries. The signal indicates the likely direction of the
root $x^*$ from $x$. If $x^* > x$, then $Z(x) = 1$ with
probability $p$ and $Z(x) = -1$ with probability $q = 1 - p$, for some fixed $p \in (1/2, 1]$.
Similarly, if $x^* < x$, then $Z(x) = -1$ with
probability $p$, and $Z(x) = 1$ with probability $q$. If $x = x^*$
then we arbitrarily assume $Z(x^*) = 1$ with probability $p$.

The PBA begins with a prior density $f_0$ on $[0, 1]$ that is positive
everywhere. Let $F_0$ denote the corresponding cumulative distribution
function. Then, for $n = 0, 1, 2, \ldots$, it iterates as follows.
\begin{enumerate}
\item Determine the next measurement point $X_n$ as the median of
  $f_n$, $X_n = F_n^{-1}(1/2)$.
\item Query at the point $X_n$ to obtain $Z_n = Z_n(X_n)$.
\item Update the density $f_n$:
\begin{align}
& \text{if $Z_n(X_n)=+1$, then }
f_{n+1}(y) = \Bigg\{
\begin{array}{ll}
2 p f_n(y),  & \; \; \mbox{if } y \geq X_n, \\
2 q f_n(y),  & \; \; \mbox{if } y < X_n,
\end{array} \label{eq:PBAUpdating1} \\
&\text{if $Z_n(X_n)=-1$, then }
f_{n+1}(y) = \Bigg\{
\begin{array}{ll}
2 q f_n(y),  & \; \; \mbox{if } y \geq X_n, \\
2 p f_n(y),  & \; \; \mbox{if } y < X_n.
\end{array} \label{eq:PBAUpdating2}
\end{align}
\item $n \leftarrow n+1$
\end{enumerate}
The updating of the density is quite natural, scaling up the mass in
the perceived direction of the root, and scaling down the mass in the
other direction. In fact, this update gives the Bayesian posterior
distribution of the location of the root under appropriate
assumptions; see \cite{waefrahen13}.

The assumption that the signal has constant probability of being
correct, irrespective of how close $x$ is to $x^*$, is problematic in
root-finding problems, and we now weaken that assumption.

\section{The Algorithm.}\label{sec:algorithm}
Now suppose that when we query $g$ at a point $x \in [0, 1]$, we
obtain $Y(x) = g(x) + \epsilon(x)$, where
$\epsilon(x)$ is a random variable, independent of any previous
function evaluations, that represents the observation noise. Repeated
observations at a fixed point $x$ yield i.i.d.\ observations $(Y_i(x)
= g(x) + \epsilon_i(x): i = 1, 2, \ldots, m)$.

Let $\tilde Z(x) = 2 \one(Y(x) \ge 0) - 1$, so that $\tilde Z(x) = 1$
if $Y(x) \ge 0$ (so that we believe $x^*$ to be to the right of $x$)
and $\tilde Z(x) = -1$ if $Y(x) < 0$. Also, define
$$\tilde p(x) = \begin{cases}
P(\tilde Z(x) = 1) & \text{if } x \le x^*,\\
P(\tilde Z(x) = -1) & \text{if } x > x^*,
\end{cases}$$
so that $\tilde p(x)$ gives the probability that $\tilde Z(x)$ correctly
indicates the direction of the root $x^*$ from $x$.

\begin{assumption}\label{assume:p bigger}
For any $x \ne x^*$, $\tilde p(x) > 1/2$.
\end{assumption}

Assumption~\ref{assume:p bigger} can be assured through conditions on
the distribution of the noise random variables $\epsilon(x)$. For
example, the assumption holds if $\epsilon(x)$ has, for all $x$, median 0 and a strictly
positive density component in a neighborhood of 0, as would be the case with unbiased
normally distributed estimators of $g(x)$. To see why we might need
the positive density condition or a similar regularity condition,
suppose that $\epsilon(x)$ has, for all $x$, a density that is
positive only on $(-\infty, -1) \cup (1, \infty)$ with (non-unique)
median 0. In this case, if $|g(x)| < 1$ then $p(x) = 1/2$.

We would like to use $\tilde p(x)$ in place of $p$ in the update step of the PBA (\ref{eq:PBAUpdating1}--\ref{eq:PBAUpdating2}), but
$\tilde p(x)$ is unknown. Furthermore, under reasonable assumptions about the
errors $\epsilon(x)$, when $g$ is continuous at $x^*$ one expects
that $\tilde p(x) \to 1/2$ as $x \to x^*$, so that the density updates 
become less and less effective as $x \to x^*$.

To overcome this difficulty, we replace the signal $\tilde Z(x)$ with
a new signal, $Z(x)$, which has probability $p(x)$ of being
correct. Moreover, $p(x) \ge p_c > 1/2$ for all $x \in [0, 1]
\backslash \{x^*\}$ for some $p_c$ chosen by the user, so that $p(x)$
is bounded away from $1/2$. Even though $p(x)$ is unknown, $p_c$ {\em
  is} known, and we can then update the density using
\refeq{eq:PBAUpdating1}--\refeq{eq:PBAUpdating2} using some $p \in
(1/2, p_c)$. This
update does not represent a true Bayesian update to a posterior
distribution, but can be viewed as a natural extension of the PBA for
problems with nonconstant and unknown $p(\cdot)$.

To obtain $Z(x)$ for $x \ne x^*$, we repeatedly evaluate the function
$g$ at $x$, observing a
sequence of signs $(\tilde Z_i(x))_i$. If $x < x^*$ then $E[\tilde Z_i(x)]=2
\tilde{p}(x) - 1 > 0$ and if $x > x^*$ then $E[\tilde Z_i(x)]=1 - 2
\tilde{p}(x) < 0$. Our goal at a particular $x \ne x^*$ is then to
evaluate the sign of the drift of the simple random walk $S(x) =
(S_m(x): m \ge 0)$ with $S_0(x) = 0$ and
$S_m(x)=\sum_{i=1}^m \tilde Z_i(x)$ for $m \ge 1$.  Sequential tests of
power one provide a powerful tool to decide whether the drift $\theta$
of a random walk satisfies the hypothesis $\theta<0$ versus
$\theta>0$. 

In a slight abuse of notation, let $S(\theta)$ denote a generic random
walk with increments in $\{-1, 1\}$ and drift $\theta$. A test of
power one can be defined
through a positive sequence $(k_m)_m$ and a stopping time
$\cN(\theta)=\inf \{ m \in \natch: |S_m(\theta)| \geq k_m \}$
\cite{Waeber2013}.  The test
decides that $\theta > 0$ if \mbox{$S_{\cN(\theta)}(\theta)
  \geq k_{\cN(\theta)}$}, that $\theta < 0$ if
$S_{\cN(\theta)}(\theta) \leq -k_{\cN(\theta)}$ and does not make a
decision if $\cN(\theta)=\infty$.  Furthermore, for a chosen
confidence parameter $\gamma \in (0,1)$ such a test satisfies
$P(\cN(\theta) < \infty) \leq \gamma$ if $\theta=0$ and
$P(\cN(\theta)< \infty)=1$ if $\theta \neq 0$.
In Section~\ref{sec:proofs} we show that one may take
\begin{equation}\label{eq:k_m}
k_m = (2m(\ln(m+1) - \ln \gamma))^{1/2}
\end{equation}
This is not an optimal test. Indeed, as shown more generally in
\cite{RoSi74,La77}, the expected value of the stopping time of our
test when the true drift is $\theta \ne 0$ is $O(\theta^{-2}
\ln(|\theta|^{-1}))$ as $\theta \to 0$, whereas \cite{Fa64} establishes
that an optimal test should have expected stopping time of the order
$O(\theta^{-2} |\ln|\ln|\theta|||)$. Unfortunately, optimal tests are
difficult to implement in practice because they rely on constants that are difficult to
compute. The power-one test we choose is sufficient for
our purposes.

One might expect that a more powerful signal could be obtained by
working directly with the observations $Y(x) = g(x) + \epsilon(x)$,
rather than with the sign of the observations. One could envisage
using other power-one tests built on repeated iid observations of
$Y(x)$, but we expect \cite{Fa64,La77,RoSi74} that while this might
accelerate the method, it would not
change the order of the rate of convergence. Accordingly, we
work with the construction above.

Assume now that the PBA measures at some point $X_n \neq x^*$ at the
$(n+1)$st iteration.  We then observe a random walk with $m$th term
$S_{n,m}=S_{n,m}(X_n)=\sum_{i=1}^m Z_{n,i}(X_n)$ until the test of
power one terminates.  Denote the stopping time of the power one test
by $N_n$, which is almost surely finite since $X_n \neq x^*$, and
define the new signal
\begin{equation*}
{Z}_n(X_n) =
\Bigg\{
\begin{array}{ll}
{}+1, & \; \; \mbox{if } S_{n,N_n} > 0, \\
{}-1, & \; \; \mbox{if } S_{n,N_n} < 0.
\end{array}
\end{equation*}

Recall that $S(0)$ is a zero-drift simple random walk. On the event $X_n > x^*$,
\begin{align}
\nonumber P(Z_n(X_n) = + 1 | X_n) & = P ( S_{n,N_n} > 0, N_n <
                                      \infty | X_n)  \\
\nonumber & \leq P( S_{\cN(0)}(0) > 0, \cN(0) < \infty) \\
\label{eq:signalpositive} & \leq \gamma/2,
\end{align}
where the first inequality follows by a sample path argument and the
second inequality by the property that $P(\cN(0)<\infty) \leq \gamma$.  An
analogous argument shows that on the event $X_n < x^*$,
\begin{align}
\label{eq:signalnegative}
P \left({Z}_n(X_n) = -1 \; \middle\vert \; X_n \right) \leq \gamma/2.
\end{align}
So, for $x \in [0,1] \setminus \{ x^* \}$, we have ${p}(x) \geq
1-\gamma/2$, where
\begin{align*}
{p}(x) =
\begin{cases}
P \left( {Z}_n(x) = +1 \right), & \; \; \mbox{if } x < x^*, \\
P \left( {Z}_n(x) = -1 \right), & \; \; \mbox{if } x > x^*.
\end{cases}
\end{align*}
Defining $p_c=1-\gamma/2$ it follows that ${p}(x) \geq p_c$ for all $x \neq x^*$, where $p_c \in (1/2,1)$ is a chosen constant (since one can choose $\gamma \in (0,1)$ in the construction of the test of power one).

It remains to discuss the case where $X_n = x^*$.  In this case the
test of power one might not terminate and the search algorithm might
stall.  From a theoretical point of view this is not convenient
because the sequence $(X_n)_n$ is then not well-defined for all $n \ge
0$.  In practice, the stalling of the algorithm can be desirable since
in this case the point $x^*$ is successfully located. (Of course,
stalling at $x^*$ cannot be detected, because it is impossible to test
in finite time whether the
stopping time is finite or not.) The event that $X_n=x^*$ for some
finite $n$ seems very unlikely in practice. Consider, for example, the
case where $x^*$ is a realization of a random variable with a positive
density on $[0,1]$:  Then the probability that $X_n=x^*$ for any $n
\in \natch$ is zero, since $X_n$ can only assume values on a set of
cardinality $2^n$.  We therefore adopt Assumption~\ref{ass:neq}
below. This can be assured, e.g., by artificially inflating the
interval $[0, 1]$ to $[0, 1 + \varepsilon]$ for some $\varepsilon$ that
is uniformly distributed on $[0, \zeta]$ for any $\zeta > 0$, 
taking $\tilde p(x) = 1$ for $x > 1$, and finally rescaling so that
the interval $[0, 1+\varepsilon]$ is mapped back to $[0, 1]$.

\begin{assumption} \label{ass:neq} $X_n \neq x^*$ for all~$n = 0, 1, 2, \ldots$ almost surely.
\end{assumption}

To summarize the algorithm then, we begin with a uniform prior on the
location of the root, and then repeatedly perform a power-one test at
the median, updating the prior using some $p \in (1/2, p_c)$, where
$p_c$ is the lower bound on $p(\cdot)$ used in the power-one test.

Using the PBA with signals $({Z}_n(X_n))_n$ instead of $(\tilde
Z_n(X_n))_n$ introduces two time scales, namely, the \emph{macro time}
$n$ corresponding to the number of different measurement points
$(X_n)_n$ and the \emph{wall-clock time} $T$ counting the total number
of function evaluations.

Intuitively, the closer the current measurement point $X_n$ is to
$x^*$, the closer $\tilde p(X_n)$ is to 1/2, and the longer the test
of power one requires to terminate.  In fact, the expected hitting
time is larger than $O((\tilde p(x) - 1/2)^{-2})$ \cite{La77,RoSi74}.  So, if $\tilde p(x) \to 1/2$ as $x \to x^*$, as would
typically hold if $g$ is continuous at $x^*$, then the number of
function evaluations between two macro iterations is likely to become
very large, which will lead to a discrepancy in the convergence
behaviors in the two time scales.


\section{Geometric Convergence.}\label{sec:exponential}
Recall that $X_n$ gives the $n$th measurement point, with $X_n$ being
the median of the distribution $\mu_n$ associated with the $n$th
density $f_n$, $n \ge 0$. In this section
we establish geometric a.s.\ convergence of $X_n \to x^*$. This geometric
convergence arises in the macro time scale, i.e., the quantity $n$
gives the number of power-one tests that are applied.

\begin{theorem}\label{th:main}
Suppose Assumptions 1-3 hold and
that we begin with a uniform prior. Then there exists $r > 0$ such
that $e^{rn}(X_n - x^*) \to 0$ as $n \to \infty$ a.s.
\end{theorem}

\proof{Proof of Theorem~\ref{th:main}.}
We assume that $x^* \in (0, 1)$. A very
similar, but simpler, argument that we omit can be employed to prove
the result if $x^* \in \{0, 1\}$. Fix $a \in (0, x^*)$ and define $A = [0, x^* - a], B = (x^* - a, x^*], C =
(x^*, 1]$. Define $T(a) = \inf\{n' \ge 0: X_n \in B \cup C \quad \forall n \ge n'\}$ to
be the time required until the sequence of medians never re-enters
$A$. Let $\stochle$ denote stochastic ordering, so that $X \stochle Y$
iff $P(X>x) \le P(Y>x)$ for all $x$. We will show that
\begin{equation}\label{eq:the whole thing}
T(a) \stochle \tau(a) + R,
\end{equation}
where $\tau(a)$ and $R$ are independent random variables with
appropriately light tails, and the nonnegative random
variable $R$ does not depend on $a$. This implies that for some $r > 0$, $e^{rn} |X_n - x^*|
I(X_n \le x^*) \to 0$ as $n \to \infty$ a.s., since for arbitrary
$\epsilon > 0$,
\begin{align}
P(e^{rn}|X_n - x^*| \one(X_n \le x^*) > \epsilon) &= P(x^* - X_n >
\epsilon e^{-rn}) \nonumber\\
&\le P(T(\epsilon e^{-rn}) > n) \nonumber\\
& \le P(\tau(\epsilon e^{-rn}) + R > n) \nonumber\\
& \le P(\tau(\epsilon e^{-rn}) > n/2) + P(R > n/2). \label{eq:log plus constant}
\end{align}
Lemma~\ref{lem:R light tailed} shows that $R$ has an exponentially
decaying tail and thus $P(R > n/2)$ is finitely summable in
$n$. Lemma~\ref{lem:B supermg} establishes that the first term in
\eqref{eq:log plus constant} is finitely summable in $n$. The
Borel-Cantelli lemma allows us to conclude that $e^{rn} |X_n - x^*| \one(X_n \le x^*) \to
0$ as $n \to \infty$ a.s. The same overall argument can be used to give
the symmetric result that $e^{rn} |X_n - x^*| \one(X_n > x^*) \to
0$ as $n \to \infty$ a.s.$\Halmos$
\endproof

The remainder of this section is dedicated to the supporting results
needed to establish Theorem~\ref{th:main}. The constant $a$ and sets
$A, B$ and $C$ defined in the proof of Theorem~\ref{th:main} are
fixed and employed throughout the remainder of this section. In the
argument, we
decompose the time till we never return to $A$ by a sum of two
terms. The first term, $\tau(a)$, gives the time required to get
sufficient conditional mass in the set $B$, i.e., sufficient mass in
$B$ relative to the mass in $A \cup B$. From this point, we show that 
there will be a finite number of excursions during which the
conditional mass in $A$ climbs back up above $1/2$, before the
conditional mass in $A$ never climbs that high again, and hence
the sequence of medians never returns to $A$.

Let $(\Omega, \cF, P)$ be an underlying probability space with a filtration $\cG
= (\cG_n: n \ge 0)$.
Let $\mu_n$ denote the random measure associated with the $n$th
density $f_n$, $n \ge 0$. So for each fixed $n \ge 0$ and $\omega \in
\Omega$, $\mu_n(\cdot, \omega)$ is a probability measure on the Borel
sets of $[0, 1]$, and for each Borel set $D \subset [0, 1]$, $\mu_n(D,
\cdot)$ is a $\cG_n$-measurable random variable. The median $X_n$ of
$\mu_n$ is $\cG_n$ measurable, and is uniquely defined because $f_n$
is strictly positive, for all $n \ge 0$. The response $Z_n(X_n)$ is
$\cG_{n+1}$ measurable for each $n \ge 0$.

We first show that the sequence of medians spends linear time on the
left side of the root $x^*$ with very high probability. Fix $\delta
\in (0, 1/2)$ and, for $n \ge 1$, let
$$N(n) = \sum_{k=0}^{n-1} \one(\mu_k([0, x^*]) \ge 1/2 + \delta)$$
be the number of instances from time 0 to time $n-1$ when the mass of $A
\cup B$ is at least $1/2 + \delta$. (It is important that the upper
index $n-1$ is not $n$ when we define certain supermartingales below.) Define
$N(0) = 0$.

\begin{lemma}\label{lem:N is linear}
We can choose $\delta > 0$ in the definition of $N(\cdot)$ above such
that there exist constants $\alpha > 0, r_1 > 0$ 
depending on $x^*$ and $\delta$ but not on $a$ such that $P(N(n) \le
\alpha n) \le \varphi(\mu_0(A \cup B)) e^{-r_1 n}$. The function
$\varphi(\cdot)$ is given in \refeq{eq:varphi def}, does not
depend on $a$, and in particular is bounded above when $\mu_0(A \cup B)$ is
bounded away from 0.
\end{lemma}

Define $\nu_n$ to be the (random) measure corresponding to the
conditional posterior distribution, conditional on lying to the left
of $x^*$, so that for any measurable $D$, $\nu_n(D) = \mu_n(D \cap [0,
x^*]) / \mu_n([0, x^*])$. Next, fix $\Delta \in (0, 1/2)$ 
and define the stopping time $\tau = \tau(a)$ to be
the first time that the conditional mass in $B$ lies above $1-\Delta$
and the median $X_n$ lies to the left of $x^*$, i.e.,
$$\tau = \inf\{n \ge 0: \nu_n(B) > 1-\Delta, X_n < x^*\}.$$
For real numbers $b, c$, let $b \wedge c = \min\{b, c\}$.

\newcommand{\sigman}{\mathcal{G}_{n}}
\begin{lemma}\label{lem:supermg1}
For $n \ge 0$, define $M_n = e^{rN(n \wedge \tau)} / \nu_{n \wedge
  \tau}(B)$. Suppose that $p(\cdot) \ge p_c > p > 1/2$. Then there exists $r = r_2 > 0$ such that
$(M_n: n \ge 0)$ is a supermartingale with respect to the filtration
$\cG$. The constant $r_2$ depends only on $\Delta, \delta, p_c$
and $p$.
\end{lemma}


\begin{lemma}\label{lem:B supermg}
Under the uniform prior, $E e^{r_2 N(\tau)} \le x^* / a$ and 
the first term in \eqref{eq:log plus constant} is summable.
\end{lemma}

\proof{Proof. }
From the supermartingale property and the uniform prior, we have that
$$E e^{r_2 N(n \wedge \tau)} \le E \frac{e^{r_2 N(n \wedge
    \tau)}}{\nu_{n \wedge \tau}(B)} = E M_n \le E M_0 = \frac 1 {\nu_0(B)}
  = \frac{x^*}{a}.$$
But $N(\cdot)$ is nondecreasing, so $e^{r_2N(n \wedge \tau)} \nearrow
e^{r_2 N(\tau)}$ as $n \to \infty$. Monotone convergence then yields
\begin{equation}\label{eq:B supermg ineq}
E e^{r_2 N(\tau)} = \lim_{n \to \infty} E e^{r_2 N(n \wedge \tau)}
\le x^*/a.
\end{equation}

For the particular choice $a = \epsilon e^{-rn}$, we obtain $E e^{r_2 N(\tau(\epsilon e^{-rn}))} \le e^{rn} x^*/\epsilon$, and hence 
\begin{align}
P(\tau(\epsilon e^{-rn}) > n/2) & \le P(N(\tau(\epsilon e^{-rn})) \ge
N(n/2)) \label{eq:Ntrans}\\
& \le P(N(\tau(\epsilon e^{-rn})) \ge N(n/2), N(n/2) > \alpha n/2) +
P(N(n/2) \le \alpha n/2) \nonumber \\
& \le P(N(\tau(\epsilon e^{-rn})) > \alpha n/2) + \beta e^{-r_1 n/2}
\nonumber\\
& \le e^{-r_2 \alpha n/2} Ee^{r_2 N(\tau(\epsilon e^{-rn}))} + \beta
e^{-r_1 n/2} \label{eq:marky}\\
& \le e^{-r_2 \alpha n/2} \frac{x^*}{\epsilon}e^{rn} + \beta
e^{-r_1 n/2}.\nonumber
\end{align}
Here we used the increasing nature of $N(\cdot)$ in \refeq{eq:Ntrans}
and a Markov bound in \refeq{eq:marky}. Taking $r < \alpha r_2/2$ gives the desired summability.
$\Halmos$\endproof

We now need to show the summability of the second term in
\eqref{eq:log plus constant}. Define the stopping time $\eta$ to be
the first time that the conditional mass in $A$ reaches 1/2 or more, i.e.,
$\eta = \inf\{n \ge 0: \nu_n(A) \ge 1/2\}$. Implicit in this
definition is the fact that $\eta = \infty$ if $\nu_n(A) < 1/2$ for
all $n$. (The specific value $1/2$ is chosen so that for $n < \eta$
we have $\mu_n(A) \le \nu_n(A) < 1/2$, so that $X_n \in B \cup C$ for
$n < \eta$.)

\begin{lemma}\label{lem:supermg2}
For $n \ge 0$, define $L_n = e^{rN(n \wedge \eta)} \nu_{n \wedge
  \eta}(A)$. There exists $r = r_3 > 0$ such that
$(L_n: n \ge 0)$ is a supermartingale with respect to the 
filtration $\cG$. The constant $r_3$ depends only on $p_c, p$ and $\delta$.
\end{lemma}

Define $U_0 = \tau = \inf\{n \ge 0: \nu_n(B) \ge
1-\Delta, X_n < x^*\}$. Lemma~\ref{lem:B supermg} establishes that
$U_0$ is finite a.s. Now, for $n \ge 1$, recursively define
\begin{align*}
V_n &= \inf\{k > U_{n-1}: \nu_k(A) \ge 1/2\}, \text{ and}\\
U_n &= \inf\{k > V_n: \nu_k(B) \ge 1-\Delta\}.
\end{align*}
Here $V_n$ is the first time after time $U_{n-1}$ that the conditional
mass in $A$ becomes at least 1/2, and $U_n$ is the first time after
$V_n$ that the conditional mass in $B$ is once again at least
$1-\Delta$. Implicit in the definition of $U_n$ and $V_n$ for $n \ge 1$ is that
they take the value $\infty$ if the stated event does not occur. Let
$\Gamma$ be the number of ``cycles,'' i.e., $\Gamma =
\sum_{n=1}^\infty \one(V_n < \infty)$.

We can now write
$$T(a) = U_0 + \sum_{j=1}^\Gamma [(V_j - U_{j-1}) + (U_j - V_j)].$$
In this expression, $V_j - U_{j-1}$ represents the $j$th time required to
increase the conditional mass in $A$ from below $\Delta$ to 1/2 or more, while $U_j - V_j$ is the subsequent time required to
ensure that the conditional mass in $B$ is at least $\Delta$. We take the empty sum (when $\Gamma = 0$) to equal 0.

Consider $U_j - V_j$, conditional on $\Gamma \ge j$, i.e.,
conditional on both $U_j$ and $V_j$ being finite. This quantity
represents the number of steps required to return the conditional mass
in $B$ to at least $1-\Delta$, starting from a point where the
conditional mass in $A$ is at least 1/2.

\begin{lemma}\label{lem:smalljumps1}
Conditional on $\Gamma \ge n$, $U_j - V_j \stochle \xi_j$ for $1 \le j
\le n$, where
$\xi_1, \xi_2, \ldots$ are iid random variables with
exponential tails, the decay rate of which does not depend on $a$.
\end{lemma}

Our next result establishes that $\eta$ has a geometrically decaying tail.

\begin{lemma} \label{lem:prob bd}
For $n \ge 0$,
$$P(\eta = n) \le 2 \nu_0(A) e^{-r_3 \alpha n} + \beta e^{-r_1 n}.$$
Here the constant $\beta$ depends on $\mu_0([0, x^*])$.
\end{lemma}

\proof{Proof. }
In this result we do not assume a uniform prior, because we will apply
the result with varying initial conditions. First note that $P(\eta =
0) = \one(\nu_0(A) \ge 1/2) \le 2 \nu_0(A)$, so the result is
immediate. For $n \ge 1$,
\begin{align}
P(\eta = n) &= P(\nu_n(A) \ge 1/2, \eta \ge n) \nonumber\\
&= P(\nu_{n \wedge \eta}(A) \ge 1/2, \eta \ge n) \nonumber\\
& \le P(e^{r_3 N(n\wedge \eta)}\nu_{n \wedge \eta}(A) \ge \frac 12
e^{r_3 N(n\wedge \eta)}, N(n \wedge \eta) \ge \alpha n, \eta \ge n)
\nonumber \\
& \qquad +
P(N(n \wedge \eta) < \alpha n, \eta \ge n) \nonumber\\
& \le  P(L_n \ge \frac 12 e^{r_3 \alpha n}) + P(N(n) < \alpha n) \nonumber\\
& \le 2 E(L_n) e^{-r_3 \alpha n} + \beta e^{-r_1 n} \label{eq:mark}\\
& \le 2 E(L_0) e^{-r_3 \alpha n} + \beta e^{-r_1 n} \nonumber\\
& = 2 \nu_0(A) e^{-r_3 \alpha n} + \beta e^{-r_1 n}. \nonumber
\end{align}
We used Markov's inequality in \refeq{eq:mark}.
$\Halmos$\endproof

Next consider $V_{j+1} - U_j$ on the event $V_{j+1} < \infty$, for
$j \ge 0$. This
quantity gives the number of steps required for the conditional
mass in $A$ to re-enter $[1/2, 1]$ after going below $\Delta$.

\begin{lemma}\label{lem:hitB}
Suppose $\nu_0(B) \ge 1 - \Delta$ and $\mu_0(A \cup B)$ is bounded
away from 0. Then there exists $r_4 > 0$ such
that
$$E [e^{r_4 (V_{j+1} - U_j)}; V_{j+1} < \infty]$$
is bounded for all $j \ge 0$.
\end{lemma}

\proof{Proof. }
Suppose that we re-initiate the process at time 0 under the posterior at time
$U_j$. From Lemma~\ref{lem:prob bd},
\begin{align*}
E [e^{r (V_{j+1} - U_j)}; V_{j+1} < \infty] &= \sum_{n=0}^\infty
e^{rn}P(\eta = n) \\
& \le \sum_{n=0}^\infty e^{rn} [2 \nu_0(A) e^{-r_3 \alpha n} +
\varphi(\mu_0(A \cup B))
e^{-r_1 n}].
\end{align*}
But $\nu_0(A) \le \Delta$ and $\varphi(\mu_0(A \cup B))$ is bounded,
since $\mu_0(A \cup B)$ is bounded away from 0 by assumption. Choosing $r= r_4 < \min\{r_1,
r_3\}$ gives the result.
$\Halmos$\endproof

The assumptions of Lemma~\ref{lem:hitB} are satisfied at each of the
times $U_0, U_1, \ldots$ that are finite. This follows by definition
of $U_0$, and for $j \ge 1$, it follows since the conditional mass in
$B$ can only strictly increase when the median is to the left of
$x^*$, so that the mass in $A \cup B$ immediately before
time $U_j$ is at least $1/2$ and therefore at least $q$ at time $U_j$.

Our next result establishes that when $\nu_0(A) < 1/2$, the
probability that the sequence of medians ever returns to $A$ is
bounded away from 1, so that the number of cycles $\Gamma$ is
stochastically bounded
by a geometric random variable that does not depend on $a$.

\begin{lemma}\label{lem:geom}
If $\mu_0(A \cup B) \ge 1/2$ and $\nu_0(B) > 1-\Delta$ then $P(\eta < \infty) \le \rho$ for some
constant $\rho < 1$.
\end{lemma}

\proof{Proof. }
From Lemma~\ref{lem:prob bd},
$$P(\eta < \infty) = \sum_{n=0}^\infty P(\eta = n) \le
P(\eta < K) + \sum_{n=K}^\infty [2 \nu_0(A) e^{-r_3 \alpha n} +
\beta(\mu_0(A \cup B)) e^{-r_1 n}]$$
for any fixed $K > 0$. Now, $\nu_0(A) \le \Delta$ and
  $\beta(\mu_0(A \cup B))$ is bounded because of the assumption that $\mu_0(A \cup B) \ge 1/2$. Choose
$K$ so that the infinite sum is at most $1/2$. Then choose $\Delta$
small enough that $P(\eta < K) = 0$. (This is possible since $\eta =
\inf\{n\ge 0: \nu_n(A) \ge 1/2\}$, and $\nu_n(A)$ increases by at most
a factor $p/q$ at each step, and will therefore hold if $\Delta
(p/q)^K < 1/2$.$\Halmos$
\endproof

\begin{lemma} \label{lem:R light tailed}
Under the assumptions of Theorem~\ref{th:main}, the second term in
\refeq{eq:the whole thing} has a finite moment generating function in a
neighborhood of 0.
\end{lemma}

\proof{Proof. }
Recall that
$$T(a) = U_0 + \sum_{j=1}^\Gamma [(V_j - U_{j-1}) + (U_j - V_j)].$$
Lemma~\ref{lem:geom} establishes that $\Gamma$ is stochastically
dominated by a geometric random variable, $\bar \Gamma$ say, such that
$P(\bar \Gamma = k) = (1-\rho) \rho^k$ for $k \ge 0$. Also,
Lemmas~\ref{lem:smalljumps1} and \ref{lem:hitB} together establish
that the summands are stochastically bounded by an iid sequence with
exponentially decaying tails where the decay rate does not depend on $a$. Such a geometric sum has moments of all orders which do not
depend on $a$. The geometric sum plays the role of the random variable
$R$ in \refeq{eq:the whole thing}.$\Halmos$
\endproof

\section{Near-Canonical Convergence.}\label{sec:canonical}
Theorem~\ref{th:main} establishes that the sequence $(X_n: n \ge 0)$
converges to the true root $x^*$ at an exponential rate in terms of
$n$, the number of macro replications. Each macro replication requires
running a test of power one, consisting of observing a random
walk until the random walk hits a boundary. Let $N_i$ denote the
number of steps of the $i$th random walk associated with macro
replication~$i$, $i \ge 0$. We complete the macro replication
associated with $X_n$ by time $T_n = \sum_{i=0}^n N_i$, $n \ge 0$.

The quantity $T_n$ is a better measure of computational complexity
than $n$, accounting as it does for the computation in each macro
replication. We show in this section that when we
measure time using $T_n$, which we refer to as ``wall-clock
time,'' we should typically expect a convergence
rate that is slower than, but arbitrarily close to, that of stochastic
approximation.

We first show that we should expect the convergence rate to remain exponential in
wall-clock time when $\tilde g(x)$ is bounded away from 0 for $x \ne
x^*$, as in certain edge-detection problems \cite{CaNo08a} and in
certain non-smooth simulation-optimization problems. In the latter
case, \cite{Li11} establishes that under modest conditions a stochastic approximation
algorithm can be expected to produce estimates of $x^*$ with far
larger error that is on the order $O(T^{-1})$, where $T$ is the total
computational time.

\begin{proposition}\label{prop:discts exp}
In addition to the assumptions of Theorem~\ref{th:main}, suppose that
$\tilde p(x)$ is bounded away from $1/2$ for $x \ne x^*$,
i.e., there exists $p_c > 1/2$ such that $\tilde p(x) \ge p_c$ for all
$x \ne x^*$. Then there exists $r > 0$ such that $e^{r T_n} |X_n -
x^*| \to 0$ as $n \to \infty$ a.s.
\end{proposition}

\proof{Proof.} The result follows from the
observation that for all $i \ge 0$, conditional on $X_i$, $N_i \stochle N$, where $N$ is
a random variable depending only on the lower bound $p_c$, with
appropriately light tails. Suppose $X_i=x < x^*$. (The case $X_i > x^*$
is similar.) Then the random walk $S = (S_k: k \ge 0)$ associated with the power-one test
has positive drift $2 \tilde p(x) - 1 \ge 2 p_c - 1 > 0$. Then the
stopping time $N_i$ of the test \refeq{eq:power one test} is stochastically bounded by the
hitting time, $N$, to the (positive only) boundary $(k_m: m \ge 1)$ of the random
walk $S$, which has increments equal to $1$ with probability $p_c$ and $-1$ with
probability $1-p_c$.

The random variable $N$ has a finite moment generating function in a
neighbourhood of zero. Indeed, for $m \ge 1$,
\begin{align*}
P(N > m) & \le P(S_m \le k_m) \\
&= P((p-q) - S_m/m \ge (p-q) - k_m / m) \\
&\le \exp[-m(p-q-k_m/m)^2/2],
\end{align*}
where the final inequality uses Hoeffding's inequality, and applies
provided that $p-q-k_m/m > 0$. But $k_m$ is of order $(m \ln
m)^{1/2}$, so for sufficiently large $m$, $P(N > m)$ decreases
exponentially rapidly, and in particular, $N$ has finite mean.

Now, from Theorem~\ref{th:main} there is some $r_1 >
0$ such that $e^{r_1 n} |X_n - x^*| \to 0$ as $n \to \infty$
a.s. Thus, it suffices to show that there is some $r > 0$ such that
$e^{rT_n}/e^{r_1 n} \to 0$ as $n \to \infty$ a.s., i.e., that $r T_n -
r_1 n \to -\infty$ as $n \to \infty$. But this again follows from
stochastic domination, since $T_n$ is stochastically dominated by the
sum of $n+1$ i.i.d.\ random variables, each of which is distributed
according to $N$, and we may then choose $r < r_1 / E(N)$.
$\Halmos$
\endproof

The situation is more difficult if $\tilde p(x)$ approaches $1/2$ as
$x \to x^*$, which is probably typical. For example, this happens if
$g$ is continuous and the error distribution at any $x$ is normal with
mean 0 and constant variance. In the case where $\tilde p(x)$
approaches 1/2 sufficiently quickly we have the following negative result, based on a result
for tests of power one from \cite{Fa64}, showing that we should
expect the error $X_n - x^*$ to be stochastically larger than $T_n^{-1/2}$, which would be the
order of the error if $X_n$ converged at the canonical rate as typically achieved
by stochastic approximation.

\begin{proposition}\label{prop:NotTight}
In addition to Assumptions 1--3, suppose that there exists $\epsilon, c > 0$
such that $0 \le \tilde p(x) - 1/2 \le c|x - x^*|$ for all $x \in (x^* - \epsilon, x^*
+ \epsilon)$. Let $(X_n: n \ge 0)$ be the sequence of medians, where
we use a test of power one to determine the updating step. Then $(T_n^{1/2}|X_n - x^*|: n \ge 0)$ is not tight.
\end{proposition}

\proof{Proof.}
It suffices to show that $(T_n(X_n - x^*)^2: n \ge 0)$ is not
tight. Since $T_n \ge n$, the result is immediate if $X_n \not \to x^*$ a.s.\
as $n \to \infty$. So suppose $X_n \to x^*$ as $n \to \infty$ a.s.
The stated condition on $\tilde p(x)$ can be extended to all $x \in [0,
1]$ through modification of the constant $c$, and so for all
$n$,
\begin{align*}
T_n(X_n-x^*)^2 & \ge N_n (X_n - x^*)^2 \\
 & \ge N_n \frac{(\tilde p(X_n)-1/2)^2}{c^2}.
\end{align*}
Recall that $N_n$ is the time required for the test of power one to
conclude. That test is based on a random walk with drift $2\tilde p(X_n)
- 1$, so the result follows if we show that $(N(\theta_n) \theta_n^2:
n \ge 0)$ is not tight, where $N(\theta)$ is the stopping time in a
power-one test based on a random walk with drift $\theta$, and
$(\theta_n: n \ge 0)$ is an arbitrary positive sequence where
$\theta_n\to 0$ as $n \to \infty$.

By passing to a subsequence, we may assume that $\rho_1 k_n / n \le \theta_n \le \rho_2
k_n/n$, where $k_n$ is the value of the
upper boundary of the continuation region of the power-one test at
time $n$ and $\rho_1, \rho_2 \in (0, 1)$ are constants that enable us to employ a result
from \cite{Fa64} below. Indeed, if we denote $D_>$ to be the event that we
decide that $\theta > 0$ in the power-one test, then from Lemma~4 of
\cite{Fa64} we have
$$\lim\inf_{n \to \infty} P(N(\theta_n) \ge n)
\ge \lim\inf_{n \to \infty} P(\{N(\theta_n) \ge n\} \cap D_>) \ge 1 - \gamma,$$
where $\gamma$ is the confidence parameter of the power-one test as
discussed in \S\ref{sec:pba}. It follows that $P(\theta_n^2
N(\theta_n) \ge n \theta_n^2) \ge 1 - 2\gamma$ infinitely often.
We chose $k_n$ to be of order $(2n \ln n)^{1/2}$, so $n \theta_n^2$ is asymptotically of order
$\ln n$, and in particular grows without bound. Thus, for any given $\epsilon < 1
- 2 \gamma$ we will not be able to find a
finite $M = M(\epsilon)$ such that $P(\theta_n^2 N(\theta_n) \ge M) <
\epsilon$ for all $n$, which concludes the proof.
$\Halmos$
\endproof

While probabilistic bisection cannot achieve the canonical rate that
one hopes for with root-finding algorithms, as we see next it can come
arbitrarily close. To see how, observe that the runlength $N_n$
becomes large when evaluating at a point $X_n$ that is ``close'' to
the root. In effect, probabilistic bisection spends a great deal of
wall-clock time at points close to the root. This motivates an
averaged root estimator of the form
\begin{equation}\label{eq:perfectaverage}
\frac 1 {T_{n-1}}\sum_{i=0}^{n-1} N_i X_i.
\end{equation}

In simulation experiments reported in Chapter 4 of \cite{Waeber2013} this estimator
performs strongly, but we have been unable to establish any theoretical
results on its asymptotic convergence behaviour. However, we have been
able to establish rates of convergence for a closely related
estimator. For
$\epsilon \in (0, 1/2)$, define
$$\hat X_n(\epsilon) = \frac{\sum_{i=0}^{n-1} N_i^{1/2-\epsilon}
  X_i}{\sum_{i=0}^{n-1} N_i^{1/2 - \epsilon}}.$$
The estimator $\hat X_n(\epsilon)$ places large weight on points $X_i$
at which the runlength $N_i$ is large, suggesting close proximity to
the root, but the weights are less extreme than those of
\refeq{eq:perfectaverage}. We next show that this estimator has an error
of order at most $T_{n-1}^{-(1/2-\epsilon)}$. Since $T_{n-1}$ is a reasonable
proxy for the wall-clock time needed to compute the estimator, we thus
establish a near-canonical convergence rate when $\epsilon$ is small,
and since $\epsilon \in (0,1/2)$ was arbitrary, we obtain the desired
conclusion that the probabilistic bisection algorithm converges at a
near-canonical rate.

\begin{theorem}\label{th:AveragedEstimator}
Suppose Assumptions 1--3 hold,
that we begin with a uniform prior, and that $\tilde p(x)-1/2 \ge c|x-x^*|$
for some $c > 0$. Then $\{|\hat X_n(\epsilon) - x^*| T_{n-1}^{1/2 -
  \epsilon}: n \ge 0\}$ is tight.
\end{theorem}

\begin{remark}
The assumption that $\tilde p$ increases at least linearly near the
root is reminiscent of similar conditions, such as a non-zero
derivative at the root, for root-finding algorithms. Without such a
condition, $\tilde p$ could grow so slowly near the root that the
runlengths in the tests of power one could be arbitrarily large.
\end{remark}

Our proof of Theorem~\ref{th:AveragedEstimator} requires a lemma.

\begin{lemma}\label{lem:finiteas}
Let $W = (W_n: n \ge 0)$ be a sequence of non-negative random variables,
and let $(\cF_n: n \ge 0)$ be a filtration. If $\sum_{n=0}^\infty
E(W_n|\cF_n) < \infty$ a.s., then $\sum_{n=0}^\infty W_n < \infty$ a.s.
\end{lemma}

\proof{Proof.} Suppose if possible that $P(A) > 0$ where $A$ is the
event $\{\sum_{n=0}^\infty W_n = \infty\}$. Since $\sum_{n=0}^\infty
E(W_n|\cF_n) < \infty$ a.s., there exists $M > 0$ such that $P(B(W,\infty)) \ge
1-P(A) / 2$, where we define the events $B(W, m)$ for $m = 0, 1, 2, \ldots, \infty$
by
$$B(W, m) = \left\{\sum_{n=0}^m E(W_n|\cF_n) \le M\right\}.$$
Then $P(A \cap B(W, \infty)) \ge P(A) + P(B(W, \infty)) - 1 > 0$. If $I_B$ is the
indicator random variable on the event $B$, then
\begin{equation}\label{eq:inf}
E \left(\sum_{n=0}^\infty W_n I_{B(W,\infty)}\right) \ge E
\left(\sum_{n=0}^\infty W_n I_{A \cap B(W,\infty)}\right) =
\infty.
\end{equation}

On the other hand, by the monotone convergence theorem,
$$E \left(\sum_{n=0}^\infty W_n I_{B(W,\infty)}(\cdot)\right) = \lim_{m \to \infty} E
\left(\sum_{n=0}^m W_n I_{B(W,\infty)}(\cdot)\right) \le E \left(\sum_{n=0}^m W_n
I_{B(W,m)}\right),$$
where the last inequality follows since $B(W,\infty) \subseteq B(W,m)$ for
any $m<\infty$. In the remainder of the proof, we show that the right-hand
side is bounded by $M$, contradicting \refeq{eq:inf}.

We proceed by induction on $m$. The event $B(W,0)$ is $\cF_0$ measurable, so
$E[W_0 I_{B(W,0)}] = E[E(W_0|\cF_0) I_{B(W,0)}] \le E[M I_{B(W,0)}] \le M$,
so the result holds for $m = 0$.

Now suppose it holds for (finite) $m-1 \ge 0$.
The event $B(W,m)$ is $\cF_m$ measurable, so that
$$E \left(\sum_{n=0}^m W_n I_{B(W,m)}\right) = E \left(\sum_{n=0}^m
  E(W_n|\cF_m)I_{B(W,m)} \right),$$
implying that we may assume without loss of generality that $W_m$ is
$\cF_m$ measurable. Now, define $W' = (W'_0, W'_1, \ldots)$
where $W'_n = W_n$ for $n \ne m$ and
$$W_m' = \left[M - \sum_{n=0}^{m-1} E(W_n|\cF_n)\right]^+,$$
with $[x]^+ = \max\{x, 0\}$. Then
\begin{align}
E\left(\sum_{n=0}^m W_n I_{B(W,m)}\right)
& \le E\left(\sum_{n=0}^m W'_n I_{B(W,m)}\right) \label{eq:dominate on B}\\
& \le E\left(\sum_{n=0}^m W'_n I_{B(W',m)}\right) \label{eq:set contained}.
\end{align}
Here, \refeq{eq:dominate on B} follows because $W_n' = W_n$ for all $n
= 0, 1, \ldots, m-1$, and because $W_m$ is $\cF_m$ measurable so that
$W_m = E(W_m|\cF_m)$, implying that on the event $B(W,m)$, $W_m' \ge
W_m$. The inequality \refeq{eq:set contained} follows since $B(W,m)
\subseteq B(W',m)$.

Now define $W'' = (W''_0, W''_1, \ldots)$ by $W''_n = W'_n$
for $n \ne m-1$ and $W''_{m-1} = W'_{m-1} +
W'_m$. Then
$$E\left(\sum_{n=0}^m W'_n I_{B(W',m)}\right) = E\left(\sum_{n=0}^{m-1} W''_n I_{B(W'',m-1)}\right),$$
which, by the inductive hypothesis, is bounded by $M$. $\Halmos$
\endproof

\proof{Proof of Theorem~\ref{th:AveragedEstimator}.}
Let $\delta = 1/2-\epsilon \in (0, 1/2)$. Then
\begin{align}
T_{n-1}^\delta |\hat X_n(\epsilon) - x^*|
&= \frac{T_{n-1}^\delta}{\sum_{i=0}^{n-1} N_i^\delta} \left|
  \sum_{i=0}^{n-1} N_i^\delta (X_i - x^*)\right| \nonumber\\
& \le \frac{T_{n-1}^\delta}{\left(\sum_{i=0}^{n-1} N_i\right)^\delta} 
  \sum_{i=0}^{n-1} N_i^\delta |X_i - x^*| \nonumber\\
&= \sum_{i=0}^{n-1} N_i^\delta |X_i - x^*| \nonumber\\
&\le \sum_{i=0}^\infty N_i^\delta |X_i - x^*|.\label{eq:DeltaBound}
\end{align}
We will show, using Lemma~\ref{lem:finiteas}, that \refeq{eq:DeltaBound}
is a finite-valued random variable a.s.,, thereby establishing
that $T_n^\delta |X_n - x^*|$ is uniformly dominated (in $n$),
completing the proof.

First, recall that $N_i$ is the sample size used in the power-one test
at $X_i$, the conditional distribution of which, conditional on $X_i$, depends only on the drift $2
\tilde p(X_i) - 1$. Let $N(\theta)$ generically denote the sample
size of the power-one test when the drift is $\theta$. Then \cite{RoSi74,La77},
$$\limsup_{\theta \to 0} E N(\theta)\theta^2 (\ln|\theta^{-1}|)^{-1} <
\infty.$$
Moreover, a sample-path argument establishes that $E N(\theta)$ is
decreasing in $|\theta| > 0$, and so it follows that for any $\gamma >
0$ there is a $\theta_0 = \theta_0(\gamma) > 0$ such that for $\theta
\ne 0$,
\begin{equation}\label{eq:ShortTest}
E N(\theta) \le c_1 |\theta|^{-(2+\gamma)} I(|\theta|
\le \theta_0) + c_2 I(|\theta| \ge \theta_0).
\end{equation}
Here and in what follows, $c_1, c_2, \cdots$ denote constants that,
apart from their existence, do not bear on the argument.

At the point $X_i$ with success probability $\tilde p(X_i)$,
$\theta = \theta(X_i) = 2 \tilde p(X_i) - 1 \ge 2c|X_i -x^*|$
where the inequality is by assumption. Recalling that $X_i$ is $\cG_i$
measurable, it follows from \refeq{eq:ShortTest} that for $i \ge 0$,
$$E[N_i|\cG_i] \le c_3 |X_i-x^*|^{-(2+\gamma)} I(|X_i - x^*| \le c_4) +
c_5 I(|X_i - x^*| > c_4).$$
Jensen's inequality then allows us to conclude that
\begin{equation}\label{eq:ENbound}
E[N_i^{\delta}|\cG_i] \le c_3 |X_i-x^*|^{-\delta(2+\gamma)} I(|X_i - x^*| \le c_4) +
c_6 I(|X_i - x^*| > c_4).
\end{equation}

Let $r > 0$ be such that $e^{rn}|X_n - x^*| \to 0$ as $n \to \infty$
a.s., as in Theorem~\ref{th:main}. Next define the index sets $J(1) =
\{i \ge 0: |X_i - x^*| > c_4\}$, $J(2) = \{i \ge 0: i \notin J(1),
|X_i - x^*| > e^{-ri}\}$ and $J(3) = \{i \ge 0: i \notin J(1) \cup
J(2)\}$.

Choose $\gamma > 0$ so that $(1-\delta(2+\gamma)) > 0$.
Using \refeq{eq:ENbound}, we obtain
\begin{align*}
E[N_i^{\delta}|X_i-x^*| \big| \cG_i] & \le I(i \in J(1)) c_6 |X_i -
x^*| +  I(i \in J(2)) c_3
  |X_i-x^*|^{1 - \delta(2+\gamma)} \\
& \qquad \quad +  I(i \in J(3)) c_3 e^{-r(1 -
                                       \delta(2+\gamma))i}\\
& \le I(i \in J(1)) c_6 +  I(i \in J(2)) c_3 +  I(i \in J(3)) c_3 e^{-r(1 - \delta(2+\gamma))i},
\end{align*}
since $|X_i-x^*| \le 1$.
Summing over $i \ge 0$ and noting that $r(1 - \delta(2 + \gamma)) >
0$, we see that
\begin{equation}\label{eq:LastStep}
\sum_{i =0}^\infty E[N_i^{\delta}|X_i-x^*| \big| \cG_i] \le c_6
|J(1)| + c_3 |J(2)| + c_7.
\end{equation}
Theorem~\ref{th:main} shows that $J(1)$ and $J(2)$ are finite sets a.s., and
the result then follows from \refeq{eq:LastStep},
\refeq{eq:DeltaBound} and Lemma~\ref{lem:finiteas}.
\endproof

\section{Additional Proofs} \label{sec:proofs}

\proof{A Power One Test}
Here we establish that \refeq{eq:k_m} has the property that when the
drift of the simple random walk $(S_m: m \ge 0)$ with $S_0 = 0$ is
$\theta = 0$, then the stopping time is finite with probability at
most $\gamma$. To that end, Robbins \cite{Ro70} showed via a likelihood ratio
test he attributes to Ville \cite{Vi39} that if $B_m = \sum_{i=1^m} Z_i$,
with $(Z_i: i \ge 1)$ being iid Bernoulli random variables with $P(Z_i
= 1) = p = 1 - P(Z_i =0)$, then $P(N'(p) < \infty) \le \gamma$, where
$$N'(p) = \inf\left\{m \ge 1: {m\choose B_m} p^{B_m} (1-p)^{m-B_m} \le
  \gamma / (m+1)\right\}.$$

The ``stopping'' region for a given $m$ and Bernoulli increments is thus the set of values $k$ such
that $P(B_m = k)$ is less than $\gamma / (m+1)$. This set is difficult to explicitly
characterize. But for $k \ge mp$, $P(B_m = k) \le P(B_m \ge k)$, and
we can apply Hoeffding's inequality to obtain an upper bound on the
latter. This yields a sufficient condition for $P(B_m = k) \le
\gamma / (m+1)$. Thus we can define a new stopping time $N''(p)$ such
that $N''(p) \ge N'(p)$, and then $P(N'(p) < \infty) \le 
P(N''(p) < \infty) \le \gamma / (m+1)$. To that end, Hoeffding's
inequality yields, for $k \ge mp$,
$P(B_m \ge k) \le \exp(-2m(k/m - p)^2)$
which is bounded above by $\gamma / (m+1)$ provided that
$k \ge mp + (m (\ln(m+1) - \ln(\gamma))/2)^{1/2}$.
Similarly, for $k \le mp$, Hoeffding's inequality applied to the left
tail yields $P(B_m = k) \le P(B_m \le k) \le \exp(-2m(k/m - p)^2)$,
which yields a symmetric result. Hence we can take
$$N''(p) = \inf\{m
\ge 1: |B_m - mp| \ge  (m (\ln(m+1) - \ln(\gamma))/2)^{1/2}\}.$$

We can now translate this particular stopping time to the desired
stopping time for the random walk $(S_m: m \ge 0)$ with zero drift
$\theta$, because $S_m = 2 B_m - m = 2(B_m - mp)$ with $p =
1/2$. Accordingly, we can take as stopping time
\begin{equation}\label{eq:power one test}
\cN(0) = \inf\{m \ge 1: |S_m| \ge  k_m\},
\end{equation}
where $k_m = (2m (\ln(m+1) - \ln(\gamma)))^{1/2}$.
\endproof

\proof{Proof of Lemma~\ref{lem:N is linear}.}
The condition $\mu_n([0, x^*]) \ge 1/2 + \delta$ is equivalent to
$\mu_n([x^*, 1]) \le 1/2 - \delta$ (since $\mu_n(\{x^*\}) = 0$ for all
$n$), which in turn is equivalent to $Y_n \ge -\ln(1/2-\delta)$, where
$$Y_n = \begin{cases}
\ln \mu_n([0, x^*]) & \text{if } \mu_n([0, x^*]) < 1/2 \\
-\ln \mu_n((x^*, 1]) & \text{if } \mu_n([0, x^*]) > 1/2 \\
0 & \text{if } \mu_n([0, x^*]) = 1/2.
\end{cases}
$$

The stochastic process $(Y_n: n \ge 0)$ takes negative values when the
median is to the right of $x^*$, positive values when the median is to
the left of $x^*$ and equals 0 when the median equals $x^*$. Assuming
that with probability 1 $X_n \ne x^*$ $\forall n$, we can ignore the
knife-edge case where $Y_n = 0$.

Consider the case where $Y_n < 0$. Then on the next step, conditional
on $X_n$ and $Y_n$, the mass in
$[0, x^*]$ is multiplied by $2p$ with probability $p(X_n)$ and by $2q$ with
probability $q(X_n) = 1 - p(X_n)$. Let $\bar p$ and $\bar q$ be shorthand
for $p(X_n)$ and $q(X_n)$ respectively. Hence, if $2p \mu_n[0, x^*]
< 1/2$, i.e., if $\ln(2p) + Y_n \le \ln(1/2)$ then $X_n$ will remain
in $[x^*, 1]$ on the next step and so, conditional on $X_n$ and $Y_n$,
\begin{equation}\label{eq:negative increm}
Y_{n+1} - Y_n = \begin{cases}
\ln(2p) & \text{w.p. }\bar p, \\
\ln(2q) & \text{w.p. }\bar q.
\end{cases}
\end{equation}
However, if $Y_n < 0$ but $\ln(2p) + Y_n \ge \ln(1/2)$, then $X_n$
will not necessarily remain in $[x^*, 1]$ on the next step. In this
case, $Y_n$ can change sign, so that conditional on $X_n$ and $Y_n$,
\begin{equation}\label{eq:cross1}
Y_{n+1} - Y_n =  \begin{cases}
-\ln(1 - 2pe^{Y_n}) - Y_n & \text{w.p. }\bar p, \\
\ln(2q) & \text{w.p. }\bar q.
\end{cases}
\end{equation}

Similarly, the increment distribution for the case when $Y_n > 0$
conditional on $X_n$ and $Y_n$ is, for $Y_n - \ln(2p) \ge  -\ln(1/2)$,
\begin{equation}\label{eq:positive increm}
Y_{n+1} - Y_n = \begin{cases}
-\ln(2p) & \text{w.p. }\bar p, \\
-\ln(2q) & \text{w.p. }\bar q,
\end{cases}
\end{equation}
and if $Y_n - \ln(2p) < -\ln(1/2)$ is
\begin{equation}\label{eq:cross2}
Y_{n+1} - Y_n = \begin{cases}
\ln(1 - 2pe^{Y_n}) - Y_n & \text{w.p. }\bar p, \\
-\ln(2q) & \text{w.p. }\bar q.
\end{cases}
\end{equation}

We will construct a new stochastic process $(\tilde Y_n: n \ge
0)$ such that $\tilde Y_n \le Y_n$ for all $n$, and such that
$\tilde N(n) = | \{k=0, \ldots, n-1: \tilde Y_n \ge
-\ln(1/2-\delta)\}|$ satisfies
$$P(\tilde N(n) \le \alpha n) \le \varphi(\mu_0(A \cup B)) e^{-r_1 n}$$
for all $n$. This will complete the proof, because $N(n) \ge \tilde
N(n)$ by pathwise domination, and then
$P(N(n) \le \alpha n) \le P(\tilde N(n) \le \alpha n)$. We construct
this process in two steps.

In the first step we replace $\bar p$ with $p$ and $\bar q$ with $q$ in both \refeq{eq:negative
  increm} and \refeq{eq:cross1}. Next, replace $\bar p$ with 1 and $\bar q$
with 0 in both \refeq{eq:positive increm} and
\refeq{eq:cross2}. Finally, replace the increment $-\ln(2p)$ in
\refeq{eq:positive increm} (which has probability 1) with $-(Y_n +
\ln(1/2))$. This last change ensures that if the process starts above
level $-\ln(1/2) + \ln(2p)$, it will immediately transition down to
$-\ln(1/2)$.

After implementing these changes, the new process is a
Markov chain on $(-\infty, \infty)$ that is pathwise dominated
by $Y$, provided that the new process is initiated at time 0 below
$Y_0$. The chain also has positive drift on the interval $(-\infty,
-\ln(1/2) - \ln(2p))$. This chain is difficult to analyze, however, because it is not
$\phi$-irreducible. Accordingly, in the second step we modify the initial value of
the chain along with all of the values of the increments by rounding
them down to nearby rational values, so as to ensure that the final
chain, $\tilde Y = (\tilde Y_n: n \ge 0)$ lives on the rationals. To
this end, let $m$ be a large positive integer that we will specify shortly, and
set $\tilde Y_0 = \lfloor m Y_0\rfloor / m$. Then $\tilde Y_0 \le Y_0$
and is rational. Furthermore, replace all of the values in the
increment distributions with similar values. For example, in
\refeq{eq:positive increm} we replace $\ln(2p)$ with $\lfloor m
\ln(2p) \rfloor / m$ and alter $\ln(2q)$ similarly. Once this is done
for all four increment distributions we obtain a process $\tilde Y$
that is stochastically dominated by $Y$, i.e., $\tilde Y_k \le Y_k$
for all $k \ge 0$. If we now choose $m$ large enough, then $\tilde Y$
has positive drift on the interval $(-\infty, \ln(1/2) - \ln(2p))$.

At this point, we have obtained a Markov chain $\tilde Y$ that lives
on a subset of the rationals and is dominated by $Y$. Moreover it is $V-$uniformly ergodic, where $V(x) = e^{|x|/2}$ for $x < 0$, $V(x) = 1$ for $0 \le x < -\ln(1/2) + \ln(2p)$
and $V(x) = \sqrt{2}(p^{1/2} + q^{1/2})^{-1}$ for $x \ge -\ln(1/2) +
\ln(2p)$. To see why, note that if $Y_0 = y < \ln(1/2) - \ln(2p)$,
then
\begin{align*}
E V(Y_1) & = pe^{(|y| - \ln(2p))/2} + q e^{(|y| - \ln(2q))/2} \\
& = V(Y_0) (p(2p)^{-1/2} + q(2q)^{-1/2}) \\
&= V(Y_0) (\sqrt{p/2} + \sqrt{q/2}).
\end{align*}
The factor $\zeta = \sqrt{p/2} + \sqrt{q/2} < 1$ since $\zeta^2 = 1/2
+ \sqrt{pq}$ and $pq < 1/4$. Also, if $Y_0 = y > -\ln(1/2) + \ln(2p)$
then $E V(Y_1) = 1 = \zeta V(Y_0)$.

Let $\pi$ be the stationary distribution of the chain $\tilde Y$. From \cite[Theorem~6.3]{konmey03}, we 
conclude that for any subset $D$ of the real line with $\int_D \pi(dx)
\in (0, 1)$ and any $\epsilon > 0$,
$$P\left(\left|\frac 1 n \sum_{k=0}^{n-1} \one(\tilde Y_k \in D) -
    \pi(D)\right| > \epsilon\right) \le \tilde \varphi(\tilde Y_0) e^{-r_1 n},$$
for certain constants $\tilde \varphi(\tilde Y_0) > 0$ and $r_1 > 0$. If
we now take $D = [-\ln(1/2-\delta), \infty)$ for sufficiently small and positive
$\delta$, we find that $\pi(D) \in (0, 1)$ and thus
$$P(|N(n)  - n \pi(D)| > \epsilon n) \le \tilde \varphi(\tilde Y_0) e^{-r_1 n}.$$
Taking $\epsilon = \pi(D) / 2$ and defining $\alpha = \pi(D) / 2$ we
obtain the desired result, except for the conclusion about the
function $\varphi$. Now, from \cite[Theorem~4.1]{konmey03}, the quantity
$\tilde \varphi(\tilde Y_0)$ is bounded by $c_1 V(\tilde Y_0)$, which is
in turn bounded by $c_2V(Y_0)$, since $Y_0 -
\tilde Y_0 \in [0, 1/m)$. The constants $c_1$ and $c_2$
do not depend on $a$. So if $w = \mu_0(A \cup B) \ge 1/2$ then $Y_0 =
-\ln(\mu_0(C))$ and thus $\tilde Y_0 \ge 0$ so that $\tilde \varphi(\tilde
Y_0) \le \zeta^{-1}$. So for $w \ge 1/2$, we can take $\varphi(w) =
c_1 c_2\zeta^{-1}$. If $w = \mu_0(A \cup B) < 1/2$ then $Y_0 = \ln(w)$
and thus $\tilde Y_0 \in [w - 1/m, w)$. It follows that we can take
$\varphi(w) = c_1 c_2 V(w -
1/m)=c_1 c_2 \frac{e^{1/(2m)}}{w^{1/2}}$. In summary,
\begin{equation} \label{eq:varphi def}
\varphi(w) = \begin{cases}
c_1 c_2 \zeta^{-1} & w \ge 1/2 \\
c_3 w^{-1/2} & w < 1/2
\end{cases}
\end{equation}
for constants $c_1, c_2, c_3$ that do not depend on $a$.$\Halmos$
\endproof

\proof{Proof of Lemma~\ref{lem:supermg1}.}
Measurability is immediate and integrability follows since
$\nu_n(B)$ can decrease by a factor of at most $q/p$ at each step
(since $\mu_n(B)$ decreases by at most a factor $2q$ and $\mu_n(A \cup
B)$ increases by at most a factor $2p$).

For the supermartingale property itself, when $n \ge \tau$ $M_n =
M_{n+1}$ so there is nothing to show. So assume for the remainder of
the proof that $n < \tau$. If $X_n \in C$ then $M_n = M_{n+1}$ and the
supermartingale property holds. It remains to check the cases $X_n \in
A$ and $X_n \in B$.

Suppose that $X_n \in A$. Define $a_n = \mu_n([0, X_n])= 1/2$, $b_n =
\mu_n((X_n, a])$, and $c_n = \mu_n(B)$. If $Z_n = 1$, i.e., the test concludes that
the root lies to the right of $X_n$, then
$$\nu_{n+1}(B) = \frac{2pc_n}{2q a_n + 2p(b_n + c_n)} =
\frac{2pc_n}{2p(a_n+b_n+c_n) - 2(p-q)a_n} = \frac{\nu_n(B)}{1 -
  \frac{p-q}p \tilde a_n},$$
where $\tilde a_n = a_n / (a_n + b_n + c_n) = 0.5 / \mu_n(A \cup B)
\in (1/2, 1)$, and the last step follows by dividing both numerator
and denominator by $2p \mu_n(A \cup B)$. Similarly, if $Z_n = -1$,
then
$$\nu_{n+1}(B) = \frac{\nu_n(B)}{1 + \frac{p-q}{q} \tilde a_n}.$$
Setting $\bar p = p(X_n)$ and $\bar q = q(X_n)$, we obtain
\begin{align}
E\left[ \left. \frac{\nu_n(B)}{\nu_{n+1}(B)} \right| \sigman\right]
&= \bar p (1 - (1-q/p)\tilde a_n) + \bar q ( 1+ (p/q-1)\tilde a_n)
\nonumber \\
&= 1 + \tilde a_n (1 - q/p) (-\bar p + \bar q p/q). \label{eq: Acase}
\end{align}

Now, $-\bar p + \bar q p/q \le -p_c + q_c (p/q) < 0$ and thus
\refeq{eq: Acase} is bounded above by
$$1 + 1/2 (1-q/p)(-p_c + q_c p / q) < 1.$$
So if $r_2 \le -\ln[1 +  1/2 (1-q/p)(-p_c + q_c p / q)]$ then the
supermartingale property holds when $X_n \in A$.

Suppose that $X_n \in B$. Redefine $a_n = \mu_n(A)$, $c_n = \mu_n((a,
X_n])$ and $d_n = \mu_n((X_n, x^*])$. If $Z_n = 1$, then, similar to
the case where $X_n \in A$,
$$\nu_{n+1}(B) = \frac{2q c_n + 2p d_n}{2q(a_n + c_n) + 2p d_n} =
\frac{q(c_n + d_n) + (p-q) d_n}{q(a_n+c_n+d_n) + (p-q) d_n} =
\frac{\nu_n(B) + (p/q-1)\tilde d_n}{1 + (p/q-1) \tilde d_n},$$
where $\tilde d_n = d_n / (a_n + c_n + d_n) = \nu_n((X_n,
x^*])$. Also, if $Z_n = -1$ then
$$\nu_{n+1}(B) = \frac{2p c_n + 2q d_n}{2p(a_n + c_n) + 2q d_n} =
\frac{p(c_n + d_n) - (p-q) d_n}{p(a_n+c_n+d_n) - (p-q) d_n} =
\frac{\nu_n(B) - (1-q/p)\tilde d_n}{1 - (1-q/p) \tilde d_n}.$$

For the sake of notational convenience, let $x = \nu_n(B)$ and $y =
(p/q-1)\tilde d_n$. Since $\tilde d_n \le \nu_n(B)$, we have $0 \le y
\le (p/q-1)x$, and $x \le 1$. Then
\begin{align}
E\left[\left.\frac{\nu_n(B)}{\nu_{n+1}(B)}\right| \sigman\right]
&= \bar p \frac{x + xy}{x+y} + \bar q \frac{x - \frac q p xy}{x -
  \frac q p y} \nonumber\\
& = \bar p\left(x + \frac{x(1-x)}{x+y}\right) + \bar q \left(x +
  \frac{x(1-x)}{x-\frac q p y}\right) \nonumber\\
&= x + x(1-x)\left(\frac{\bar p}{x+y} + \frac{\bar q}{x - \frac q p y}\right). \label{eq:last}
\end{align}
Holding $x$ fixed, \refeq{eq:last} is strictly convex in $y$, so the maximum
value over $y \in[0, (p/q-1)x]$ is attained at one of the
endpoints of the interval. So either $y = 0$, in
which case the value of \refeq{eq:last} is 1, or $y = (p/q-1)x$ in
which case the value of \refeq{eq:last} is
$$x + (1-x)(\bar p q / p + \bar q p/q)  \le x + (1-x)(p_c q/p +
q_c p/q),$$
where the inequality arises because both sides contain convex combinations
of the same two values, with the right-hand side placing larger weight
on the larger value. We also have that
\begin{equation}\label{eq:smallerthan1}
p_c q/p + q_c p/q < p q/p + q p/q = 1,
\end{equation}
for the same reason, and so $\theta < 1$ where $\theta$ is the
constant on the left-hand side in \refeq{eq:smallerthan1}. Hence, when
$y = (p/q-1)x$, \refeq{eq:last} is at most $x + \theta (1-x) \le
1-\Delta + \theta \Delta < 1$. Of the two endpoints, the one with $y =
0$ yields the higher bound of $1$. So \refeq{eq:last} is at most equal
to 1, thereby proving the supermartingale property when $N(n+1) = N(n)$.

When $N(n+1) = N(n) + 1$, i.e., $\mu_n[0, x^*] \ge 1/2 + \delta$, then
we want to show that the quantity \refeq{eq:last} is bounded away from 1. To this end, 
$\tilde d_n \ge \delta / \mu_n([0, x^*]) \ge \delta$ so that $y \ge
(p/q-1)\delta \ge (p/q-1)\delta x$. Thus when we maximize
\refeq{eq:last}, we can now do so over the range $(p/q-1)\delta x \le
y \le (p/q-1) x$. As above, this is a maximization of a strictly
convex function over a convex set, so the maximum is attained at an
endpoint. The endpoint $y = (p/q-1)x$ yields the bound $1 - \Delta +
\theta \Delta < 1$ as before. The endpoint $y = (p/q-1)\delta x$ gives a value for
\refeq{eq:last} of
\begin{equation}\label{eq:last2}
x + (1-x)\left( \frac{\bar p}{1+(p/q-1)\delta} + \frac{\bar q}{1
    - (1-q/p)\delta}\right).
\end{equation}
If we consider this quantity as a strictly convex function of $\delta \in [0, 1]$,
the maximum is strictly attained at one of the endpoints, and hence since $\delta \in (0, 1/2)$, for some constant $\theta' \in (0,
1)$ the value of \refeq{eq:last2} is at most $x + (1-x) \theta' \le 1
- \Delta + \theta' \Delta < 1$. We have thus shown that
\refeq{eq:last} is bounded away from 1.$\Halmos$
\endproof

\proof{Proof of Lemma~\ref{lem:supermg2}.}
If $n \ge \eta$ then there is nothing to prove, so suppose that $n <
\eta$, and in particular, $\nu_n(A) < 1/2$. Then $\mu_n(A) < 1/2$ so
$X_n \notin A$. If $X_n \in C$ then $L_{n+1} = L_n$, so we need only
consider the case where $X_n \in B$. Let $a_n = \mu_n(A)$, $c_n =
\mu_n((x^*-a, X_n])$ and $d_n = \mu_n((X_n, x^*])$.

If $Z_n = +1$, then
$$\nu_{n+1}(A) = \frac{2q a_n}{2q(a_n+c_n) + 2p d_n} = \frac{a_n}{a_n
  + c_n + d_n + (p/q-1)d_n} = \frac{\nu_n(A)}{1 + (p/q-1) \tilde
  d_n},$$
where $\tilde d_n = d_n / (a_n + c_n + d_n)$ is the conditional mass
in $(X_n, x^*]$. Similarly, if $Z_n = -1$ then
$$\nu_{n+1}(A) = \frac{2p a_n}{2p(a_n+c_n) + 2q d_n} = \frac{a_n}{a_n
  + c_n + d_n - (1-q/p)d_n} = \frac{\nu_n(A)}{1 - (1-q/p) \tilde
  d_n}.$$

Hence, setting $\bar p = p(X_n)$ and $\bar q = q(X_n)$,
\begin{equation} \label{eq:easier}
E \left[ \frac{\nu_{n+1}(A)}{\nu_n(A)}\right]
= \frac{\bar p}{1 + (p/q-1) \tilde d_n} + \frac{\bar q}{1 -
  (1-q/p) \tilde d_n} = f((p/q-1)\tilde d_n),
\end{equation}
where
$$f(y) = \frac{\bar p}{1+y} + \frac{\bar q}{1-(q/p)y}.$$
Now, $f$ is strictly convex (in $y$ such that the denominators remain
positive), so it attains its maximum over any (finite) interval at an endpoint
of the interval. Also, $y = (p/q-1) \tilde d_n \le (p/q-1)/2$ since
$d_n < 1/2$ and $\tilde d_n = d_n / (1/2 + d_n)$. Thus, we need only
consider $y \in [0, (p/q-1)/2]$. Now, $f(0) = 1$ and
$$f((p/q-1)/2) = \frac{\bar p}{1/2 + p/2q} + \frac{\bar q}{1/2 +
  q/2p} = 2q \bar p + 2p \bar q \le 4pq < 1.$$
The second-to-last inequality results from the fact that $q \bar p +
p \bar q$ is a convex combination of $q$ and $p$, the size of which
can be increased by increasing the weight on the larger value, i.e.,
on $p$. Since $f$ is at most 1 at its endpoints, we have established
that \refeq{eq:easier} is at most 1.

This proves the super martingale property if $N(n+1) =
N(n)$. If $N(n+1) = N(n) + 1$, then $\tilde d_n \ge \delta$ and hence
$y$ is restricted to the interval $[(p/q-1) \delta, (p/q-1)/2]$. Over
this interval $f$ is strictly convex and thus bounded away from 1.$\Halmos$
\endproof

\proof{Proof of Lemma~\ref{lem:smalljumps1}.}
As in the proof of Lemma~\ref{eq:B supermg ineq}, we find that
$$E( e^{r_2 (N(U_n) - N(V_n))} | \mathcal{G}_{V_n}) \le 1 /
\nu_{V_n}(B).$$
At time $i = V_n - 1$, $\nu_i(B)  > 1/2$, since $\nu_i(A) < 1/2$. In a
single time step, the conditional mass in any set can be reduced by
at most a factor $q/p$, since $\mu_{i+1}(B) \ge 2q \mu_i(B)$ and
$\mu_{i+1}(A \cup B) \le 2p \mu_i(A \cup B)$. Thus
$$E( e^{r_2 (N(U_n) - N(V_n))} | \mathcal{G}_{V_n}) \le 2p/q.$$
For notational convenience, we renormalize by taking time 0 to be
the time $V_n$, and all probabilities and expectations are conditional
on $\mathcal{G}_{V_n}$. So now
\begin{align*}
P(U_n - V_n > k) &\le P(N(U_n - V_n) \ge N(k)) \\
& \le P(N(U_n - V_n) \ge N(k), N(k) \ge \alpha k) + P(N(k) < \alpha k)
\\
& \le P(N(U_n-V_n) \ge \alpha k) + \varphi(\mu_0(A \cup B)) e^{-r_1 \alpha k}\\
& \le E[e^{r_2N(U_n-V_n)}] e^{-r_2 \alpha k} + \varphi(\mu_0(A \cup B)) e^{-r_1 \alpha
  k}\\
& \le 2\frac p q e^{-r_2 \alpha k} + \varphi(\mu_0(A \cup B)) e^{-r_1 \alpha
  k}.
\end{align*}

To complete the proof we establish a uniform bound on $\varphi(\mu_0(A
\cup B))$ by showing that $\mu_0(A \cup B) \ge q$. Indeed, the
conditional distribution immediately 
prior to time $V_n$ was such that the median was in $B$, since
the conditional distribution only changes when the median is in $A
\cup B$ and the conditional mass in $A$ was less than 1/2. Thus just
prior to time $V_n$, the total mass in $A \cup B$ was at least 1/2,
this being the mass to the left of the median. Since the mass in any
set cannot decrease by more than a factor $2q$, the mass in $A \cup B$
at time $V_n$ is at least $q$.$\Halmos$
\endproof

%
%
%

\section*{Acknowledgments.}
This research was partially supported by National Science
Foundation grant CMMI-1200315.


\bibliographystyle{plain} 
\bibliography{SRFPbib} 


\end{document}